\documentclass[11pt]{amsart}

\usepackage{epsf}

\title{   Alexander duality in subdivisions of Lawrence polytopes}

\author{     Francisco Santos}
\address{    Francisco Santos,
             Departamento de Matem\'aticas, Estad\'{\i}stica y
                 Computa\-ci\'on\\
             Universidad de Cantabria, E-39005, Santander, SPAIN.}
\email{santos@matesco.unican.es}
\author{     Bernd Sturmfels}
\address{    Bernd Sturmfels,
             Department of Mathematics, University of California,
             Berkeley, CA 94720, USA.}
\email{bernd@math.berkeley.edu}

\thanks{This work was completed while the first author was visiting
  U. C. Davis, supported by U. C. Davis, M.S.R.I. and the Spanish government.
  He was also supported by grant BFM2001--1153 of the Spanish
  Direcci\'on General de Investigaci\'on. The second author was
  partially supported by NSF Grant DMS-9970254.}

\date{February 2002.}

\newtheorem{theorem}{Theorem}
\newtheorem{lemma}[theorem]{Lemma}

\newtheorem{remark}[theorem]{Remark}
\newtheorem{corollary}[theorem]{Corollary}
\newtheorem{proposition}[theorem]{Proposition}

\theoremstyle{definition}

\newcommand{\reals}{{\mathbb R}}

\newcommand{\e}{{\epsilon}}
\newcommand{\M}{{\mathcal M}}
\newcommand{\Z}{{\mathcal Z}}
\renewcommand{\H}{{\mathcal H}}

\begin{document}

\begin{abstract}
The class of simplicial complexes
representing triangulations and subdivisions of Lawrence polytopes is
closed under Alexander duality.
This gives a new geometric  model for oriented matroid duality.
\end{abstract}

\maketitle

\section{\bf Introduction}

The  aim of this note is to show that oriented
matroid duality can be seen as an instance
of Alexander duality of simplicial complexes (see e.g. \cite{BBM}).
We represent
an affine oriented matroid $(\M,f)$ on the ground set
$\{1,\ldots,n, f \}$ by a simplicial complex
$\Delta(\M,f)$ on the vertex set
$\, \{x_1,\ldots,x_n,$ $y_1,\ldots,y_n \}\,$
as follows. The facets of $\Delta(\M,f)$ are the complements
of the sets
$$  \{ x_i \, : \, i \in C^+ \} \,\, \cup \,\,
 \{ y_j \, : \, j \in C^- \}, $$
where $C = (C^+,C^-)$ runs over all signed cocircuits
of $(\M,f)$ such that the distinguished element $f$
lies in $C^+$. We have the following  result:

\begin{theorem}
\label{thm1}
The Alexander dual of $\,\Delta(\M,f)\,$ is the
simplicial complex  $\,\Delta(_{-f}\M^*,f)\,$ associated
with the affine oriented matroid $(_{-f}\M^*,f)$. Here $_{-f}\M^*$ denotes
the oriented matroid
dual to $\M$ with the element $f$ reoriented.
\end{theorem}

This duality can be expressed geometrically in terms of
Lawrence polytopes.
Suppose that the contraction
$\M/f$ is represented by a $d \times n$-matrix
${\bf D}$ of rank $d$. Then the associated Lawrence
polytope (see e.g.~\cite[\S 6.6]{Ziegler})
 is the convex hull of the columns of the
$(d + n) \times 2n$-matrix
\begin{equation}
\Lambda ( {\bf D}) \quad = \quad
\begin{pmatrix}
{\bf D} & {\bf 0} \\
{\bf I} & {\bf I}
\end{pmatrix}.
\label{eqn:lawrence}
\end{equation}
Here ${\bf I}$ is the $n \times n$-identity matrix,
${\bf 0}$ is the $d \times n$-zero matrix, and the columns
are indexed by $\, \{x_1,x_2,\ldots,x_n,y_1,y_2,\ldots,y_n \}$.
Recall that $\{x_i, y_i\}$ is the complement of a facet
of $\Lambda({\bf D})$, for all $i$.
It turns out that $\Delta(\M,f)$ is a polyhedral subdivision
of the Lawrence polytope $\Lambda({\bf D})$, where each
maximal face in the subdivision is represented  by the simplex
on its set of vertices. This subdivision is a triangulation
if and only if the matroid $\M\backslash f$ is uniform.
The Lawrence polytope $\Lambda({\bf D})$ itself is called {\em uniform}
if all $d \times d$-minors of ${\bf D}$ are nonzero, or,
in the non-realizable case, if
the matroid $\M/f$ is uniform.

The following is our main result:

\begin{theorem}
\label{thm2}
The following families of simplicial  complexes on
the $2n$-element set
$ \,\{x_1,\ldots,x_n,y_1,\ldots,y_n \}\,$
are closed under Alexander duality:
\begin{itemize}
\item [(1)] Regular triangulations of uniform Lawrence polytopes,
\item [(2)] regular subdivisions of  Lawrence polytopes,
\item [(3)] triangulations of uniform Lawrence matroid polytopes,
\item [(4)] subdivisions of Lawrence matroid polytopes.
\end{itemize}
\end{theorem}

Moreover, Alexander duality gives a bijection between
regular triangulations of Lawrence polytopes
and regular subdivisions of uniform Lawrence polytopes. These
two families are not closed under Alexander duality.

The families (3) and (4) in Theorem \ref{thm2}
 refer to the case when the oriented
matroid $\M/f$ cannot be represented by a matrix ${\bf D}$.
For the relevant definitions and notations used here
we refer to the books \cite{OMbook} and \cite{OMtri}.
In particular, see  \cite[\S 9.3]{OMbook} for Lawrence (matroid) polytopes
and \cite[\S 9.6]{OMbook} for subdivisions
of (matroid) polytopes. The first author proved in
\cite[Theorem 4.14]{OMtri} that every
subdivision of a Lawrence (matroid) polytope is induced by
a lifting of oriented matroids $\,\M/f  \longrightarrow \M$.

\begin{figure}[ht]
  \epsfxsize = 4 in
  \leavevmode
  \epsfbox{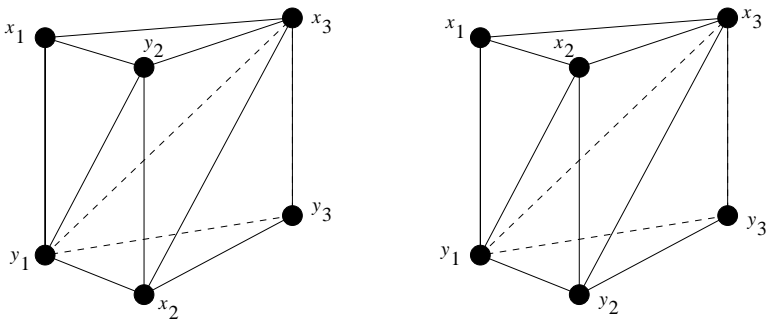}
\caption{The triangulation of a triangular prism is Alexander
  self-dual, after relabeling the vertices.
 The non-edges on the left are the complements
  of the tetrahedra on the right.}
\label{fig:prisms}
\end{figure}

Our presentation is organized as follows. In Section \ref{sec2}
we prove Theorem \ref{thm1} and we interpret
$\Delta(\M,f)$ in terms of hyperplane arrangements.
The proof of Theorem \ref{thm2} is given in Section \ref{sec4}.
Examples of Alexander dual pairs of subdivided Lawrence polytopes
are given in Section \ref{sec3}.
The smallest non-trivial example is the 
pair of triangular prisms
in Figure \ref{fig:prisms}.

Section \ref{sectopology} concerns the Alexander duals of simplicial balls and
spheres in general. This section
was added after we received the very helpful
comments of an anonymous referee. He or she pointed us to the work of Dong
\cite{Dong} and proposed the extension stated in part 2 of Theorem
\ref{thm:Dong}.

The original motivation for this project came from
commutative algebra and hyperk\"ahler geometry. The
simplicial complex $\Delta(\M,f)$ is represented
algebraically as a
squarefree monomial ideal
in $k[x_1,\ldots,x_n,y_1,\ldots,y_n]$.
The minimal free resolution of this ideal
constructed in \cite{NPS} can be interpreted
as a (suitably homogenized) coboundary complex on
the Alexander dual $\,\Delta(_{-f}\M^*,f)$. In particular,
part (1) in Theorem \ref{thm2} furnishes a large
class of Stanley-Reisner rings
which are  Cohen-Macaulay and have an
explicit linear resolution.
The quotient of such a Stanley-Reisner ring
modulo a linear system of parameters was  shown in
 \cite{hyperkaehler} to equal
the cohomology ring of a toric hyperk\"ahler variety.
These varieties are complete intersections in the
toric variety whose fan is a cone over  $\Delta(\M,f)$.
It would be interesting to explore the
duality of toric hyperk\"ahler varieties arising
from our results.

\section{Oriented Matroid Duality is Alexander Duality}
\label{sec2}

We recall the
combinatorial definition of Alexander duality.
Let $K$ be a simplicial complex on the
vertex set $V$. Then the {\em Alexander dual} of $K$
is the simplicial complex
\[
K^\vee \quad := \quad \{V\backslash \sigma : \sigma\not\in K \}
\]
The Alexander Duality Theorem states that the $i$-th reduced
homology group $\widetilde H_i(K, {\bf Z})$
of $K$ equals the $(|V|-3-i)$-th reduced cohomology group
$\widetilde H^{|V|-3-i}(K^\vee, {\bf Z})$ of $K^\vee$. See, e.g.,
\cite[equation (2)]{BBM} or \cite[(9.17)]{Bj}.
In particular, the Alexander dual of an acyclic simplicial
complex is acyclic, although
the Alexander dual of a contractible simplicial complex need not be
contractible. See Section~\ref{sectopology} for a discussion of this and
related topological issues.

\begin{proof}[Proof of Theorem \ref{thm1}]
The statement can be rephrased as the following claim:
given an oriented matroid $\M$  on the ground set
$\{1,\dots,n,f\}$, for any pair of subsets
$\sigma_1, \sigma_2\subseteq \{1,\dots,n\}$ one and only one of the
following happens:
\begin{itemize}
\item[(1)] There is a cocircuit $(C^+,C^-)$ in $\M$
   with $C^-\subseteq \sigma_1$, and $f\in C^+\subseteq \sigma_2\cup\{f\}$, or

\item[(2)] There is a cocircuit $(D^+,D^-)$ in $\M^*$
 (that is, a circuit in $\M$)
   with $f\in D^-\subseteq \{1,\dots,n,f\}\backslash \sigma_1$
and  $D^+\subseteq
   \{1,\dots,n\}\backslash \sigma_2$.
\end{itemize}
Indeed, condition (1) above is equivalent to
\[
 \{ x_i \, : \, i \not\in \sigma_2 \} \,\, \cup \,\,
 \{ y_j \, : \, j \not\in \sigma_1 \} \quad \in \quad \Delta(\M,f),
\]
and condition (2) is equivalent to
\[
 \{ x_i \, : \, i \in \sigma_2 \} \,\, \cup \,\,
 \{ y_j \, : \, j \in \sigma_1 \} \quad \in  \quad \Delta(_{-f}\M^*,f).
\]
The claim follows from Lemma \ref{lemma:duality} below, taking $e=f$
and color classes
$\, B=(\sigma_2\backslash\sigma_1)\cup\{f\}$,
$\, W=\sigma_1\backslash\sigma_2$,
$\, R=\sigma_1\cap\sigma_2$,
and $\,G=\{1,\dots,n\}\backslash(\sigma_1\cup\sigma_2)$.
We also set $\,(C^+,C^-)=(Y^+,Y^-)$ and
 $\,(D^+,D^-)=(X^-,X^+)$.

Lemma \ref{lemma:duality} is just a rephrasing
of the 4-painting axiom of oriented matroid
circuits and cocircuits. The notation in the lemma
is chosen to exactly match the axiom as it appears in
\cite[Theorem 3.4.4]{OMbook}. This is the
reason why we have $X=-D$ above rather than
reorienting $X$ in the lemma.
\end{proof}

\begin{lemma}
\label{lemma:duality}
Let $B$, $W$, $G$ and $R$ be a partition of the ground set of an
oriented matroid $\M$. Let $e\in B\cup W$ be one of the elements.
Then, exactly one of the
following happens:
\begin{enumerate}
\item[(1)] There is a circuit $(X^+,X^-)$
                    with $X^-\subseteq W\cup G$
and $e\in X^+\subseteq B\cup G$, or
\item[(2)] There is a cocircuit $(Y^+,Y^-)$
                    with $e\in Y^+\subseteq B\cup R$
and $Y^-\subseteq  W\cup R$.
\end{enumerate}
\end{lemma}

\medskip

We now interpret $\Delta(\M,f)$
in terms of hyperplane arrangements.
By the Topological Representation Theorem \cite[\S 4]{OMbook},
 an affine oriented matroid $(\M,f)$ of rank $d$
on $\{1,\dots,n,f\}$ represents an affine arrangement $\H(\M,f)$
of $n$ pseudo-hyperplanes in $\reals^{d-1}$, with the distinguished
element $f$ playing the role of the hyperplane at infinity.
We can regard $\H(\M,f)$ as a cover of $\reals^{d-1}$
by $2n$ closed half-spaces $\{x_1,\ldots,x_n,y_1,\ldots,y_n\}$,
where $x_i$ and $y_i$ label respectively the positive
and negative sides of the $i$-th oriented hyperplane.
It is straightforward to check that
a subset of these half-spaces
has a non-empty intersection in $\reals^{d-1}$
if and only if the corresponding subset
of $\{x_1,\ldots,x_n,y_1,\ldots,y_n\}$
is a simplex in
$\Delta(\M,f)$. In other words:

\begin{remark}
The simplicial complex $\Delta(\M,f)$ is the nerve
of the cover of $\reals^{d-1}$ consisting of the
$2n$ closed half-spaces in the arrangement $\H(M,f)$.
\end{remark}

The facets of  $\Delta(\M,f)$ are maximal intersecting
families of closed half-spaces. They correspond to the
vertices of the arrangement $\H(\M,f)$.
The face poset of $\H(\M,f)$ appears as a subposet
in the face poset of $\Delta(\M,f)$.
A simplex $\sigma\in\Delta(\M,f)$ is called
\emph{full} if $\sigma\cap\{x_i,y_i\}\ne \emptyset$ for all $i$.

\begin{remark}
\label{prop:arrangement}
If $\M \backslash f$ is uniform, then
the face poset of $\H(\M,f)$ is anti-isomor\-phic to the poset
of full simplices of $\Delta(\M,f)$. If $\M \backslash f$
is not uniform, then
the former is a strict subposet of the latter.
\end{remark}

This implies that the oriented matroid $\M$ can be recovered
from the simplicial complex $\Delta(\M,f)$ provided $\M$
is uniform. The same statement is not true for general
oriented matroids.
For instance, consider an arbitrary arrangement of hyperplanes
which intersect in a line, and then adjoin two parallel hyperplanes
transverse to that line. Here $\Delta(\M,f)$ consists
of two simplices of the same dimension which share a common facet,
regardless of which arrangement we started with.

\section{Lawrence Polytopes in Dimension Three, Four and Five}
\label{sec3}

In Section~\ref{sec4} we are going to prove Theorem~\ref{thm2} by translating
 Theorem \ref{thm1} into the
language of subdivisions of Lawrence (matroid) polytopes. As a preparation for
that we describe in this section all the Lawrence polytopes which exist in
dimensions up to 5, and an example of our Alexander duality result involving
two Lawrence polytopes of respective dimensions 4 and 5.

We first recall the construction of Lawrence polytopes
in oriented matroid language, and then we discuss
low-dimensional Lawrence polytopes.
Let $\M$ be an oriented matroid of rank $d$ on
$\,\{1,\dots,n\}$, and let $\M^*$ be its dual.
Let $\M^*\cup (-\M^*)$ be the oriented matroid on
$\{x_1,\dots,x_n,y_1,\dots,y_n\}$ defined by labeling the $i$-th element
of $\M^*$ as $x_i$ and extending $\M^*$ by an
element $y_i$ opposite to each $x_i$.
The dual of $\M^*\cup (-\M^*)$ is called the
\emph{Lawrence oriented matroid} (or
\emph{Lawrence polytope}, since it is a matroid
polytope) of $\M$, and denoted $\Lambda(\M)$.
It has $2n$ elements and rank $d+n$.
Lawrence (matroid)
polytopes are studied in Section 9.3 of \cite{OMbook} and in Chapter 4
of \cite{OMtri}.
For example, \cite[Lemma 4.11(ii)]{OMbook} implies that
$\Lambda(\M)$ has $n-l+2c$ facets, where $c$ is the number of cocircuits
of $\M$ and $l$ the number of coloops.

Since all the oriented matroids with $d+n\le 11$ are
realizable, all Lawrence matroid polytopes of dimension at most
$10$ are honest polytopes, that is, they can be realized
by $(d+n) \times 2n$-matrices of the form  $\Lambda({\bf D})$ as
in (\ref{eqn:lawrence}). In what follows we describe
all Lawrence polytopes of dimension $d+n-1 \leq 5 $.

\begin{figure}[ht]
  \epsfxsize = 4 in
  \leavevmode
  \epsfbox{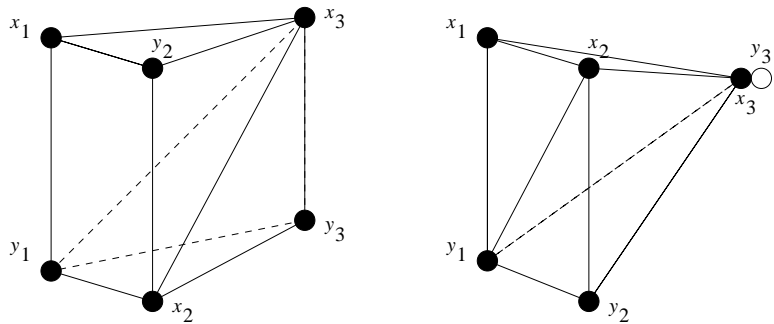}
\caption{This subdivision of a uniform Lawrence polytope
(the triangular prism) is Alexander dual to a
triangulation of a non-uniform Lawrence polytope
(the pyramid).}
\label{fig:pyramid}
\end{figure}

Let us first discuss the degenerate cases
when $\M$ has a loop or coloop.
If $x_i$ is a coloop in $\M$ (i.e.~if
the $i$-th column of ${\bf D}$ is linearly independent of all
others), then it becomes a loop in $\M^*$. Then, $x_i$ and $y_i$
are loops in $\M^*\cup (-\M^*)$ and coloops in
$\Lambda(\M)$. Geometrically,
$\Lambda({\bf D})$ is an iterated pyramid over
the Lawrence polytope $\Lambda({\bf D} \backslash \{x_i\})$.
If $x_i$ is a
loop in $\M$ (i.e.~if the $i$-th column of ${\bf D}$ is
zero), then $\Lambda(\M)$ is obtained from
$\Lambda(\M\backslash\{x_i\})$ by adjoining a pair of parallel
elements which forms a positive cocircuit. Geometrically,
$\Lambda({\bf D})$ is a pyramid over
$\Lambda({\bf D}\backslash\{x_i\})$ with
apex at a pair of identified points $x_i$ and $y_i$.
The right  picture of Figure
\ref{fig:pyramid} represents this situation.
The apex of the pyramid corresponds to the identified points
$y_3$ and $x_3$. Note that the triangulation
 uses $x_3$ and not $y_3$ as a vertex. This is
indicated in the diagram
 with a filled dot for $x_3$ and an empty dot for $y_3$.

We now consider only Lawrence polytopes that are not pyramids over
other Lawrence polytopes, which
is the same as allowing only oriented matroids without loops or coloops.
There  are eight combinatorial types of such Lawrence polytopes
having dimension at most five. The corresponding
parameters $(n,d)$ are
$\,(2,1), (3,1),(4,1),(5,1),(3,2),
 (4,2),  (4,2),  (4,2)$:

\begin{itemize}
\item If $d=1$, then the Lawrence polytope of $\M$ equals the product
$\Delta^1 \times \Delta^{n-1}$ of a segment and a simplex of dimension $n-1$.
The polytope  $\Delta^1\times \Delta^{n-1}$ has $ \,n \,! $ triangulations
each isomorphic to the well-known
staircase triangulation.
The case $n = 2$ is featured in Figure \ref{fig:prisms}.
The case $n = 3$ appears in (\ref{Staircase}) below.

\item If $n-d=1$, then $\M^*$ and $\M^*\cup(-\M^*)$ have rank 1, and
$\Lambda(\M)$ has corank 1, i.e., it has a unique circuit.
Assuming without
loss of generality that all the elements of $\M$ have the same
orientation, this unique circuit is
$(\{x_1,\dots,x_n\},\{y_1,\dots,y_n\})$. The polytope $\Lambda(\M)$
can be realized as the convex hull of the union of two
$(n-1)$-simplices in $\reals^{2n-2}$ whose
relative interiors intersect in a unique point.
This Lawrence polytope is the cyclic $(2n-2)$-polytope with $2n$ vertices.

\item
Up to reorientation, there are three oriented matroids
$\,\M_1, \, \M_2 , \,\M_3 $ of rank $2$ on $4$ elements.
They are represented by $2 \times 4$-matrices
$$ \qquad \quad
{\bf D}_1 \, = \, (v_1,v_1,v_2,v_2 ), \,\,\,
{\bf D}_2 \, = \, (v_1,v_1,v_2,v_3 ), \,\,\,
{\bf D}_3 \, = \, (v_1,v_2,v_3,v_4 ). $$
Here the $v_i$ are pairwise linearly
independent vectors in the plane.
In each case, $\Lambda({\bf D}_i)$ is a five-dimensional Lawrence
polytope
with eight vertices and with $6+2i$ facets.
For instance, $ \Lambda({\bf D}_1)$
is the join of two squares.
\end{itemize}

We shall examine the
Lawrence polytope $\,\Lambda({\bf D}_3)$ by computing
one of its triangulations along with its  Alexander dual.
We start out with the
$2 \times 5$-matrix
$$ A \quad = \quad \bordermatrix{
& x_1 & x_2 & x_3 & x_4 &  f \cr
&  1  &  1  &  1  &  1  &  1  \cr
&  4  &  3  &  2  &  1  &  0  \cr},
$$
and we fix the following  Gale dual
$3 \times 5$-matrix, with last column reoriented:
$$ B \quad = \quad \bordermatrix{
& x_1 & x_2 & x_3 & x_4 &  \overline f \cr
&  1  & -2  &  1  &  0  &  0 \cr
&  0  &  1  & -2  &  1  &  0 \cr
&  0  &  0  &  1  & -2  & -1 \cr}.
$$
Thus $A$ and $B$ represent uniform matroids.
Let $A'=A/f$ and  $B'=B/\overline f$ denote the
matrices gotten from $A$ and $B$
by contracting the last column. Contracting $f$ means
projecting every vector $v\in
A\backslash\{f\}$ along the direction of $f$ to a linear hyperplane not
containing $f$. In our case:
$$
 A' \quad = \quad \bordermatrix{
& x_1 & x_2 & x_3 & x_4 \cr
&  4  &  3  &  2  &  1  \cr},
\qquad
 B' \quad = \quad \bordermatrix{
& x_1 & x_2 & x_3 & x_4 \cr
&  1  & -2  &  1  &  0  \cr
&  0  &  1  & -2  &  1  \cr}.
$$
The $2 \times 4$-matrix $B'$ has the form of ${\bf D}_3$
in the previous paragraph and will play the role of ${\bf D}$
in the big matrix $\Lambda({\bf D})$ of equation (1).
The polytopes $\Lambda(A')$ and $\Lambda(B')$ are
$4$-dimensional and $5$-dimensional, both with eight
vertices. As we saw above, $\Lambda(A')$ is (affinely isomorphic to)
the product of a segment and a tetrahedron.

There are precisely six signed cocircuits of $B$ (or circuits of $A$)
in which the element $\overline f$ is positive:
\begin{equation}
\label{Ssix}
\{y_1, x_2, \overline f \}, \,  \{y_1, x_3, \overline f \}, \,
\{y_1, x_4, \overline f \}, \,  \{y_2, x_3, \overline f \}, \,
\{y_2, x_4, \overline f \}, \,  \{y_3, x_4, \overline f \}.
\end{equation}
There are precisely four  signed cocircuits of $A$ (or circuits of $B$)
in which the element ${f}$ is positive:
\begin{equation}
\label{FFour}
\{x_2, x_3, x_4, f \},
\{y_1, x_3, y_4, f \},
\{y_1, y_2, x_4, f \},
\{y_1, y_2, y_3, f \}.
\end{equation}
Taking complements in (\ref{Ssix}) we obtain the maximal simplices
in a  regular triangulation of the $5$-dimensional Lawrence polytope
$\Lambda(B')$:
\begin{eqnarray*}
& \{ x_1, x_2, x_3, y_1, y_2, y_4\}, \,\,
  \{ x_1, x_2, x_3, y_1, y_3, y_4\}, \,\,
  \{ x_1, x_2, x_4, y_1, y_3, y_4\}, \,\,    \\
& \{ x_1, x_2, x_3, y_2, y_3, y_4\}, \,\,
  \{ x_1, x_2, x_4, y_2, y_3, y_4\}, \,\,
  \{ x_1, x_3, x_4, y_2, y_3, y_4\}. \,\,
\end{eqnarray*}

Taking complements in (\ref{FFour}) we obtain the maximal simplices
in a staircase triangulation of the $4$-dimensional Lawrence polytope
$\,\Lambda(A') \, = \, \Delta^1 \times \Delta^3$:
\begin{equation}
\label{Staircase}
 \{ x_1, y_1, \! y_2, \! y_3, \! y_4\},
 \{ x_1, \! x_2, y_2, \! y_3, \! y_4\},
 \{ x_1, \! x_2, \! x_3, y_3, \! y_4\},
 \{ x_1, \! x_2, \! x_3, \! x_4, y_4\}.
\end{equation}
These two simplicial complexes
are Alexander dual to each other.
The Stanley-Reisner ideals
of the two triangulations are
gotten from (\ref{Ssix}) and
(\ref{FFour}) by deleting $f$ and $\overline f$ and
regarding each set as  square-free monomial. Namely,
the Stanley-Reisner ideal of our triangulation
of $\Delta(B')$ is
\begin{equation}
\label{Sssix}
\langle \,y_1  x_2, \, y_1  x_3 ,\,
 y_1 x_4 , \, y_2 x_3 , \,
y_2 x_4, \, y_3 x_4 \rangle,
\end{equation}
and  the Stanley-Reisner ideal of our triangulation
of $\Delta(A')$ is
\begin{equation}
\label{Fffour}
\langle
x_2 x_3 x_4 , \,\,
y_1 x_3 x_4 , \,\,
y_1 y_2 x_4 , \,\,
y_1 y_2 y_3 \rangle.
\end{equation}

\section{Duality of Subdivided Lawrence Polytopes}
\label{sec4}

The proof of Theorem \ref{thm2} is based on the  non-trivial fact that
all subdivisions of a Lawrence matroid polytope are
lifting subdivisions. This  fact is one of the main results  in
the monograph \cite{OMtri}.

We recall the definition of lifting subdivisions.
Let $(\M,f)$ be an affine oriented matroid on the ground set
$\{1,\dots,n,f\}$, and assume that $f$ belongs to some positive cocircuit.
Consider the sets $\, \{ \, x_i:i\not\in C^+\} \,$
where $C$ runs over all \emph{positive} cocircuits of $\M$ with  $f\in
C^+$. These sets form  (the maximal cells of) a subdivision of 
the oriented matroid $\M/f$. 
Subdivisions of an oriented matroid obtained
in this manner are called {\em lifting subdivisions}.
 For the general definition of
subdivisions of oriented matroids see
\cite[\S 9.6]{OMbook} or \cite{OMtri}.

If $\M/f$ is realized by a vector configuration, then
subdivisions of $\M/f$ are the same as polyhedral subdivisions
(also called polyhedral fans) of it. If not only $\M/f$ but also
$\M$ is realized by a vector configuration $A$,
then the lifting subdivision induced by $(\M,f)$ is
the regular subdivision of $A/f$ corresponding to the lifting
$A/f \to A$. Some lifting subdivisions of vector configurations
are not regular, and some polyhedral subdivisions are not lifting.
See \cite[Corollary 9.6.8]{OMbook}.
By \cite[Proposition 9.1.1]{OMbook},
every lifting subdivision is either a $(d-1)$-ball or 
a  $(d-1)$-sphere,
where $d$ is the rank of $\M/f$, and the latter 
happens exactly when $\M$ is acyclic and $\M/f$
totally cyclic. The topological type, or even the homotopy type,
is not known
for general subdivisions of non-realizable oriented matroids.

\begin{proposition}
\label{prop:contractible}
Let $S$ be a lifting subdivision of a rank $d$ oriented matroid on $n$
elements. If $S$ is not a triangulation we consider it as a simplicial
complex whose facets are the maximal faces of $S$.
Then, the Alexander dual $S^\vee$ of $S$ is either contractible or 
homotopy equivalent to an  
$(n-d-2)$-sphere, depending on whether $S$ itself is contractible or a
$(d-1)$-sphere.
\end{proposition}

\begin{proof}
A subset $\sigma\subseteq\{1,\dots,n\}$ is in $S^\vee$ if and only if
$\M$ has no positive cocircuit with
$f\in C^+\subseteq \sigma$. By Lemma
\ref{lemma:duality} (with
$W=R=\emptyset$, $B=\sigma$  and
$G=\{1,\dots,n,f\}\backslash \sigma$)
this happens if and only if $\M$ has a
circuit $(D^+,D^-)$ with $f\in
D^+$ and $D^-\cap \sigma=\emptyset$. Equivalently, if the
closed positive
half-spaces labeled by $\sigma$ have non-empty intersection
in the arrangement $\H(\M^*,f)$.

In other words,
$S^\vee$ is the nerve of the family of closed positive
half-spaces of $\H(\M^*,f)$. By the Nerve Theorem
(see \cite[\S 11]{Bj}) $S^\vee$ has the homotopy type of the
union of these half-spaces, which equals the complement of the (open)
cell of $\H(\M^*,f)$  corresponding to the covector
$(f,\{1,\dots,n\})$, or the entire affine
space if that covector does not appear in $\M^*$. This
complement is
contractible unless the covector exists and the corresponding cell is
bounded, in which case it is an $(n-d-2)$-sphere. The cell
$(f,\{1,\dots,n\})$ exists and is bounded if and only if
$\M^*\backslash f$ is acyclic and $\M^*$ totally cyclic.
\end{proof}

We now shift gears and replace $\M/f$ by $\Lambda (\M/f)$.
It was proved in \cite[Theorem 4.14]{OMtri} that every
subdivision of a Lawrence matroid polytope $\Lambda (\M/f)$
is a lifting subdivision. See also
\cite[\S 4]{HRS} for the realizable case. Moreover,
lifts of $\Lambda(\M/f)$ and lifts of
$\M/f$ are essentially the same thing. In particular, $(\M,f)$
represents a lift of $\Lambda(\M /f)$ and a lifting subdivision of it.
We denote this subdivision
by $S(\M,f)$. Its maximal faces are the sets
\[
\{ \, x_i:i\not\in C^+\} \, \cup \, \{\, y_i:i\not\in C^-\}
\]
where $C$ runs over all cocircuits of $\M$ with  $f\in C^+$.
Hence $S(\M,f)$ coincides with  $\Delta(\M, f)$
if we regard $S(\M,f)$ as a simplicial
complex as in the statement of Proposition \ref{prop:contractible}.
Observe that $S(\M,f)$ is a triangulation if and
only if $\M\backslash f$ is uniform.
Theorem \ref{thm1} can be rephrased as:

\begin{corollary}
\label{coro:lawrence}
Let $(\M,f)$ be an affine oriented matroid. Let $(_{-f}\M^*,f)$ be its
dual, reoriented at $f$. The subdivisions
$S(\M,f)$ and $S(_{-f}\M^*,f)$
of $\Lambda(\M/f)$ and $\Lambda(\M^*/f)$ are Alexander dual to one another.
\qed
\end{corollary}

\begin{proof}[Proof of Theorem \ref{thm2}]
Part (4) follows from Corollary \ref{coro:lawrence}. Part (3)
corresponds to the case where both $\M/f$ and $\M\backslash f$
are uniform and part (2) is the case where both $\M/f$ and
$\M\backslash f$ are realizable. Part (1) is the intersection of
both cases. Observe that $\M\backslash f$ is uniform or realizable if
and only if $\M^*/f=(\M\backslash f)^*$ has that property.
\end{proof}

Triangulations of Lawrence matroid polytopes
and subdivisions of uniform Lawrence matroid polytopes, intermediate
between cases (3) and (4) of Theorem \ref{thm2}, correspond
respectively to $\M\backslash f$ and $\M/f$ being uniform. Hence
they are not self-dual classes of simplicial complexes $\Delta(\M,f)$,
but classes dual to one another. Adding the attribute
``regular'' to both sides gives another two dual classes.
Figure \ref{fig:pyramid} was an example of this.
Figure \ref{fig:diagram} below summarizes Theorem
\ref{thm2} and this remark, showing how Alexander
duality acts on the following eight families of simplicial complexes
on $\{x_1,\dots,x_n,y_1,\dots,y_n\}$:
\begin{itemize}
\item ${\tt S} = \{$Subdivisions of matroid Lawrence polytopes$\}$.
\item ${\tt R} = \{$Regular subdivisions of Lawrence polytopes$\}$.
\item ${\tt T} = \{$Triangulations of matroid Lawrence polytopes$\}$.
\item ${\tt U} = \{$Subdivisions of uniform matroid Lawrence polytopes$\}$.
\item ${\tt RT} = {\tt R}\cap {\tt T}$,  \ \
  ${\tt RS} = {\tt R}\cap {\tt S}$, \ \
  ${\tt TU} = {\tt T}\cap {\tt U}$, \ \
  ${\tt RTU} = {\tt R}\cap {\tt T}\cap {\tt U}$.
\end{itemize}
This is a Hasse diagram: thin lines represent 
set-theoretic inclusions
among the eight families. Thick
arrows indicate the action of Alexander duality.

\begin{figure}[ht]
  \epsfxsize = 1.5 in
  \leavevmode
  \epsfbox{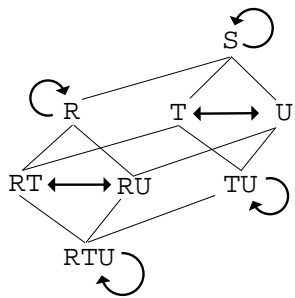}
\caption{A diagram showing the action of Alexander duality on several
  families of simplicial complexes.}
\label{fig:diagram}
\end{figure}

\begin{remark} \rm
When we say ``$\Delta(\M,f)$ is a regular
triangulation of a Lawrence polytope'' we mean
``there is a realization ${\bf D}$ of
$\M/f$ for which the subdivision corresponding to
$\Delta(\M,f)$ is regular''. A stronger meaning would be
``in every realization ${\bf D}$ of
$\M/f$ the subdivision corresponding to
$\Delta(\M,f)$ is regular''.
Theorem \ref{thm2} is not true with
this stronger meaning, as the following example shows. Let $\M$ be the
oriented matroid realized by
\[
A=\bordermatrix{
& x_1 & x_2 & x_3 & x_4 & x_5 & x_6 &  f \cr
&  1  &  2  & -\e &  0  &\e-1 & -2  &  0  \cr
& \e  &  0  &  1  &  2  & -1  & -2  &  0  \cr
&  1  &  1  &  1  &  1  &  1  &  1  &  1  \cr},
\]
where $\e$ is sufficiently small and positive.
Let $A_1=A\backslash f$ and let
$A_2=\{v_1,\dots,v_6\}$ be a realization of $\M\backslash f$ in
which the planes spanned by $\{v_1,v_2\}$, $\{v_3,v_4\}$ and $\{v_5,v_6\}$
meet in a line. Let $B_1$ and $B_2$ be Gale transforms of
$A_1$ and $A_2$, respectively.
Since $A_2$ cannot be extended to a realization
of $\M$, \  $\Delta(_{-f}\M^*,f)$ is a regular
triangulation of $\Lambda(B_1)$ but not of $\Lambda(B_2)$, even
though both represent the same matroid polytope
$\Lambda(\M^*/f)$.
On the other hand, $\Delta(\M,f)$ is a regular triangulation
of any realization of $\Lambda(\M/f)$, because any realization of
$\M/f$ is the contraction of one of $\M$.
\qed
\end{remark}

In closing we relate our discussion to {\em zonotopal tilings},
which is the geometric model  for oriented matroids
featured prominently in  \cite{Ziegler}.
Suppose that $\M/f$ can be realized as a vector
configuration ${\bf D}=\{v_1,\dots,v_n\}\subset\reals^{d-1}$.
The Bohne-Dress Theorem (see
\cite[\S 7.5]{Ziegler}) says
that the cell-complex dual to the arrangement $\H(\M,f)$
is a zonotopal tiling  $\Z(\M,f)$
of the zonotope $Z({\bf D})=\sum_{i=1}^n [O,v_i]$.
The exact relation between $\Z(\M,f)$ and $S(\M,f)$ is as follows.
Let $\pi:\Lambda({\bf D})\to \Delta^{n-1}$
be the projection sending the pair of
vertices $x_i$ and $y_i$ to the $i$-th vertex of the standard
$(n-1)$-simplex $\Delta^{n-1}$.
In coordinates, this projection just forgets the first
$d$ rows in the matrix $\Lambda({\bf D})$ given in (\ref{eqn:lawrence}).
Let $P$ be the centroid of $\Delta^{n-1}$. Then, $\pi^{-1}(P)$ is a
scaled copy of the zonotope $Z({\bf D})$. The Cayley Trick \cite{HRS}
states that the zonotopal tiling $\Z(\M,f)$ is the intersection
of the subdivision $S(\M,f)$ with that zonotope.

\section{The topology of Alexander duals.}
\label{sectopology}

We start by showing that the Alexander dual of a contractible simplicial
complex need not be contractible, with the following reasoning
suggested to us by Anders Bj\"orner.
Let $K$ be any acyclic but not contractible
simplicial complex with at least $5$ more vertices than its dimension.
Small such complexes, with dimension $2$ and $10$ vertices,
are described in \cite[p. 284]{BL}.
By the assumption on dimension,
every three vertices form a triangle in $K^\vee$, and hence 
$K^\vee$ is simply connected.
It is also acyclic by the Alexander Duality Theorem.
By standard algebraic topology results,
acyclic and simply connected simplicial complexes are contractible.

This fact contrasts the following result, pointed out to us by an
anonymous referee. Part 1 is taken from \cite{Dong}. 
The proof of the second part is due to the referee.

\begin{theorem}[Dong \cite{Dong}]
\label{thm:Dong}
Let $S$ be a simplicial complex of dimension $d$ with $n$ vertices. Then:
\begin{enumerate}
\item If $S$ is a $d$-sphere then $S^\vee$ has the homotopy 
type of the $(n-d-3)$-sphere.
\item If $S$ is a $d$-ball then $S^\vee$ is contractible.
\end{enumerate}
\end{theorem}

\begin{proof}
If $n\ge d+5$, the argument 
above gives that $S^\vee$ is simply connected. This, together with
the fact that it has the homology groups of the $(n-d-3)$-sphere
(respectively, of a contractible 
space) implies that it is homotopy equivalent to
the $(n-d-3)$-sphere (resp., it is contractible).

Let us now assume that $n\le d+4$. In part 1, this implies that $S$ is
actually polytopal, by a classical result of Mani \cite{Mani}.
Corollary 22 in \cite{Dong} implies that the
Alexander dual of a simplicial $d$-polytope with $n$ vertices
is homotopy equivalent to the $(n-d-3)$-sphere.

In part 2, the case $n\le d+3$ is proved by 
similar arguments: Coning the boundary
of $S$ to a new vertex we get a simplicial $d$-sphere with at most $d+4$
vertices, hence a polytopal one. This implies that $S$ is a shellable ball,
hence collapsible (see Lemma 17 in \cite{Dong}). The Alexander dual of a
collapsible space is contractible, by \cite[Corollary 12]{Dong}.

We still have to deal with the case $n=d+4$ in part 2. We will prove that in
this case $S^\vee$ is simply connected. Hence, the same arguments as in
the case $n\ge d+5$ apply. The complex
$S^\vee$ has a complete 1-skeleton, but not a
complete 2-skeleton. The triangles missing are precisely the complements of
the maximal simplices in $S$, and our task is to show that they all produce
null-homotopic loops. To see this, let $\sigma$ be a $d$-simplex in $S$, with
complement $\{p,q,r\}$. If
$\sigma$ has a boundary facet $\sigma\backslash\{s\}$, then $\{p,q,s\}$, 
$\{p,r,s\}$,  and $\{q,r,s\}$ are triangles in $S^\vee$, hence the loop
$\{p,q,r\}$ is null-homotopic. If $\sigma$ has no boundary facet, let
$\sigma'$ a $d$-simplex of $S^\vee$ adjacent to $\sigma$. Suppose the
complement of $\sigma'$ is $\{p,q,s\}$.
Then the triangles $\{p,r,s\}$ and $\{q,r,s\}$ are in $S^\vee$ and prove that
the loops $\{p,q,r\}$ and $\{p,q,s\}$ are homotopic. In other words, missing
triangles of $S^\vee$ corresponding to adjacent $d$-simplices of $S$ are
homotopic. Any maximal simplex in the ball $S$ can be connected to one
incident to the boundary.
This proves that every missing triangle is homotopic to
a null-homotopic one.
\end{proof}

This result in particular implies Proposition~\ref{prop:contractible} for
lifting triangulations. But actually Dong's paper~\cite{Dong} contains the
ingredients needed to generalize it to arbitrary subdivisions. Indeed, his
Theorem 27 (together with his Lemma 25)
states that the Alexander dual of
every polyhedral decomposition of a $d$-sphere, 
considered as a simplicial complex as we did in Proposition
\ref{prop:contractible}, is homotopy equivalent of a $(n-d-3)$-sphere.
But the three properties of
polyhedral complexes that he uses are also satisfied by subdivisions of
oriented matroids. Namely: (1) they are regular cell complexes, (2) the
intersection of any two closed cells is a closed cell (Dong calls this the
{\em meet
property}) and (3) they can be refined to triangulations without the addition
of new vertices by the so-called pulling construction (for the pulling
refinement of oriented matroid subdivisions see \cite[Section 9.6]{OMbook} or
\cite[Remark 4.4]{OMtri}). Hence, we can generalize Proposition
\ref{prop:contractible} as follows:

\begin{theorem}
\label{thm:contractible2}
Let $S$ be a subdivision of a rank $d$ oriented matroid on $n$
elements. If $S$ is not a triangulation we consider it as a simplicial
complex whose facets are the maximal faces of $S$.
Then, 
\begin{enumerate}
\item If $S$ (as a cell complex) is a $(d-1)$-sphere, then 
$S^\vee$ is homotopy equivalent to a  
$(n-d-2)$-sphere.
\item If $S$ (as a cell complex) is a $(d-1)$-ball, then 
$S^\vee$ is contractible.
\end{enumerate}
\end{theorem}

\begin{proof}
  Let $T$ be a triangulation obtained by pulling refinement of $S$. As
  mentioned in \cite{Dong}, $S$ (considered as a simplicial complex) collapses
  to $T$ and this implies that $T^\vee$ collapses to $S^\vee$. Since $T$ is
  homeomorphic to (the cell complex) $S$, 
 the homotopy type of $T^\vee$ is given by
  Theorem \ref{thm:Dong}.
\end{proof}

It is not known whether cases (1) and (2) of
Theorem \ref{thm:contractible2} cover all subdivisions of oriented
matroids. They cover, at least, all subdivisions of realizable ones and all
lifting subdivisions of non-realizable ones.

\end{document}